\documentclass[a4paper,11pt]{amsart}
\usepackage{graphicx}
 \usepackage{mathptmx}
\usepackage{amsmath}
\usepackage{amssymb}
\usepackage{enumitem}
\usepackage{xcolor}

\newtheorem{theorem}{Theorem}[section]

\begin{document}

\title{Some identities on type 2 degenerate Bernoulli polynomials of the second kind}

\author{Taekyun Kim $^{1}$}
\address{$^{1}$ Department of Mathematics, Kwangwoon University, Seoul 139-701, Republic of Korea}
\email{tkkim@kw.ac.kr}

\author{Lee-Chae Jang$^{2,*}$}
\address{ Graduate School of Education, Konkuk University, Seoul, 05029, Republic of Korea}
\email{Lcjang@konkuk.ac.kr}

\author{Dae San Kim$^{3}$}
\address{$^{2}$ Department of Mathematics, Sogang University, Seoul 121-742, Republic
    of Korea}
\email{dskim@sogang.ac.kr}

\author{Han Young Kim$^{4}$}
\address{$^{4}$ Department of Mathematics, Kwangwoon University, Seoul 139-701, Republic of Korea,}
\email{gksdud213@kw.ac.kr}

\maketitle

\begin{abstract}
In recent years, many mathematicians studied various degenerate versions of some special polynomials of which quite a few interesting results were discovered. 
In this paper, we introduce the type 2 degenerate Bernoulli polynomials of the second kind and their higher-order analogues, and study some identities and expressions for these polynomials. Specifically, we obtain a relation between the type 2 degenerate Bernoulli polynomials of the second and the degenerate Bernoulli polynomials of the second, an identity involving higher-order analogues of those polynomials and the degenerate Stirling numbers of second kind, and an expression of higher-order analogues of those polynomials in terms of the higher-order type 2 degenerate Bernoulli polynomials and the degenerate Stirling numbers of the first kind.

\end{abstract}

\markboth{\centerline{\scriptsize Some identities on type 2 degenerate Bernoulli 
polynomials of the second kind }}{\centerline{\scriptsize T. Kim, L.-C. Jang, D. S. Kim, H. Y. Kim}}

\section{Introduction}\label{sec1}

\vspace{1cm}

As is known, the type $2$ Bernoulli polynomials are defined by the generating function
\begin{equation}\label{eq1}
\frac{t}{e^t -e^{-t} } e^{xt} = \sum_{n=0}^{\infty}B_n^*(x)\frac{t^n}{n!},\quad \text{ (see \cite{ref7})}.
\end{equation}
From \eqref{eq1}, we note that
\begin{equation}\label{eq2}
B_n^* (x) =2^{n-1} B_n\left( \frac{x+1}{2} \right),\;\; (n\geq 0),
\end{equation}
where  $B_{n}(x)$ are the ordinary Bernoulli polynomials given by
\begin{equation*}
\frac{t}{e^t -1}e^{xt} = \sum_{n=0}^{\infty}B_n(x) \frac{t^n}{n!}.
\end{equation*}

Also, the type $2$ Euler polynomials are given by 
\begin{equation}\label{eq3}
e^{xt} \text{sech} t = \frac{2}{e^t +e^{-t} } e^{xt}
=\sum_{n=0}^{\infty}E_n^* (x) \frac{t^n}{n!}, \quad \text{ (see \cite{ref7})}.
\end{equation}
Note that
\begin{equation}\label{eq4}
E_n^* (x) = 2^n E_n \left( \frac{x+1}{2} \right),\;\;(n\geq 0),
\end{equation}
where $E_n(x)$ are the ordinary Euler polynomials given by
\begin{equation*}
\frac{2}{e^t +1} e^{xt}
=\sum_{n=0}^{\infty}E_n (x) \frac{t^n}{n!}, \quad \text{ (see \cite{ref3, ref4, ref5, 
ref8, ref10, ref11})}.
\end{equation*}

The central factorial numbers of the second kind are defined as
\begin{equation}\label{eq5}
x^n = \sum_{k=0}^n T(n,k) x^{[k]}, \;\;\;  \quad \text{ (see \cite{ref6})},
\end{equation}
or equivalently as 
\begin{equation} 
\frac{1}{k!}(e^{\frac{t}{2}}-e^{-\frac{t}{2}})^k=\sum_{n=k}^{\infty}T(n,k)\frac{t^n}{n!},
\end{equation}
where $x^{[0]}=1$, $x^{[n]}=x\left(x+\frac{n}{2}-1 \right) \left(x+\frac{n}{2}-2 \right)
\cdots \left(x-\frac{n}{2}+1 \right), \; \; (n \geq 1)$.

\vspace{0.1in}

It is well known that the Daehee polynomials are defined by
\begin{equation}\label{eq6}
\frac{\log (1+t)}{t}(1+t)^x= \sum_{k=0}^n D_n(x) \frac{t^n}{n!}, \;\;\;  \quad \text{ (see \cite{ref14, ref15})}.
\end{equation}
When $x=0$, $D_n=D_n(0)$ are called the Daehee numbers.\\

The Bernoulli polynomials of the second kind of order $r$ are defined by 
\begin{equation}\label{eq7}
\left( \frac{t}{\log (1+t)}\right)^r(1+t)^x= \sum_{k=0}^n b_n^{(r)}(x) \frac{t^n}{n!}, \;\;\;  \quad \text{ (see \cite{ref13})}.
\end{equation}
Note that $b_n^{(r)} (x)= B_n^{(n-r+1)} (x+1),\;\; (n\geq 0).$ Here $B_n^{(r)}(x)$ are the ordinary Bernoulli polynomials of order $r$ given by 
\begin{equation}\label{eq8}
\left( \frac{t}{e^t -1}\right)^re^{xt} = \sum_{k=0}^n B_n^{(r)}(x) \frac{t^n}{n!}, \;\;\;  \quad \text{ (see \cite{ref13, ref14, ref15, ref16})}.
\end{equation}
It is known that the Stirling numbers of the second kind are defined by
\begin{equation}\label{eq9}
\frac{1}{k!}\left( e^t -1\right)^k =\sum_{n=k}^\infty S_2(n,k) \frac{t^n}{n!}, \;\;\;
\quad \text{ (see \cite{ref13})},
\end{equation}
and the Stirling numbers of the first kind by 
\begin{equation}\label{eq10}
\frac{1}{k!} \log^k (1+t) =\sum_{n=k}^\infty S_1(n,k) \frac{t^n}{n!}, \;\;\;
\quad \text{ (see \cite{ref13})}.
\end{equation}

For any nonzero $\lambda\in \mathbb{R}$, the degenerate exponential function is defined by
\begin{equation}\label{eq11}
e_\lambda^x(t) = (1+\lambda t)^{\frac{x}{\lambda}} =\sum_{n=0}^\infty (x)_{n, \lambda} \frac{t^n}{n!}, \;\;\;
\quad \text{ (see \cite{ref9})},
\end{equation}
where  $(x)_{0,\lambda}=1$, $(x)_{n,\lambda} =x(x-\lambda)\cdots (x-(n-1)\lambda), \; \; (n \geq 1)$.
\\
In particular, we let
\begin{equation}\label{eq11-1}
e_{\lambda}(t) = e_{\lambda}^{1}(t) = (1+\lambda t)^{\frac{1}{\lambda}}.
\end{equation}

In \cite{ref1,ref2}, Carlitz introduced the degenerate Bernoulli polynomials which are given by
the generating function
\begin{equation}\label{eq12}
\frac{t}{e_\lambda(t) -1} e_\lambda^x (t)  = \sum_{n=0}^\infty \beta_{n,\lambda}(x) \frac{t^n}{n!}.
\end{equation}
Also, he considered the degenerate Euler polynomials given by
\begin{equation}\label{eq13}
\frac{2}{e_\lambda(t) +1} e_\lambda^x (t)  = \sum_{n=0}^\infty \mathcal{E}_{n,\lambda}(x) \frac{t^n}{n!}, \;\;\;
\quad \text{ (see \cite{ref1, ref2})}.
\end{equation}

Recently, Kim-Kim considered the degenerate central factorial numbers of the second kind given by
\begin{equation}\label{eq14}
\frac{1}{k!} \left(e_\lambda^{\frac{1}{2}}(t) - e_\lambda^{-\frac{1}{2}}(t)  \right)^k
 =\sum_{n=k}^\infty T_\lambda(n,k) \frac{t^n}{n!}.
\end{equation}
Note that $\lim_{\lambda\rightarrow 0}T_\lambda(n,k)=T(n,k)$, (see \cite{ref12}).

\vspace{0.1in}

In this paper, we introduce the type 2 degenerate Bernoulli polynomials of the second kind and their higher-order analogues, and study some identities and expressions for these polynomials. Specifically, we obtain a relation between the type 2 degenerate Bernoulli polynomials of the second and the degenerate Bernoulli polynomials of the second, an identity involving higher-order analogues of those polynomials and the degenerate Stirling numbers of second kind, and an expression of higher-order analogues of those polynomials in terms of the higher-order type 2 degenerate Bernoulli polynomials and the degenerate Stirling numbers of the first kind.

\vspace{1cm}

\section{Type $2$ degenerate Bernoulli polynomials of the second kind }\label{sec2}

\vspace{1cm}

Let $\log_\lambda t$ be the compositional inverse of $e_\lambda(t)$ in \eqref{eq11-1}.
Then we have
\begin{equation}\label{eq15}
\log_\lambda t=\frac{1}{\lambda}\left( t^\lambda -1 \right).
\end{equation}
Note that  $\lim_{\lambda\rightarrow 0} \log_\lambda t =\log t$.
Now, we define the {\it{degenerate Daehee polynomials}} by
\begin{equation}\label{eq16}
\frac{\log_\lambda (1+t)}{t} (1+t)^x = \sum_{n=0}^\infty D_{n,\lambda}(x) \frac{t^n}{n!}.
\end{equation}
Note that  $\lim_{\lambda\rightarrow 0} D_{n, \lambda}(x)=D_n(x), \;(n\geq 0)$.
In view of \eqref{eq7}, we also consider the degenerate Bernoulli polynomials of the second kind of
order $\alpha$ given by
\begin{equation}\label{eq17}
\left( \frac{t}{\log_\lambda (1+t)} \right)^\alpha (1+t)^x = \sum_{n=0}^\infty b_{n,\lambda}^{(\alpha)}(x) \frac{t^n}{n!}.
\end{equation}
Note that  $\lim_{\lambda\rightarrow 0} b_{n, \lambda}^{(\alpha)} (x)=b_n^{(\alpha)}(x), \;(n\geq 0)$.
From \eqref{eq17}, we have
\begin{equation}\label{eq18}
\left( \frac{\lambda t}{(1+t)^{\frac{\lambda}{2}}-(1+t)^{-\frac{\lambda}{2}} } \right)^\alpha (1+t)^{x- \frac{\lambda\alpha}{2}}
= \sum_{n=0}^\infty b_{n,\lambda}^{(\alpha)}(x) \frac{t^n}{n!}.
\end{equation}
For $\alpha=r\in \mathbb{N}$, and replacing $t$ by $e^{2t}-1$ in \eqref{eq18}, we get
\begin{align}\label{eq19}
\sum_{m=0}^\infty b_{m,\lambda}^{(r)} (x) \frac{1}{m!}(e^{2t}-1)^m
=& \left(  \frac{\lambda t}{e^{t\lambda}-e^{-t\lambda}} \right)^r \frac{1}{t^r}(e^{2t}-1)^r e^{(2x-\lambda r)t}\cr
=& \sum_{k=0}^\infty B_k^*\left( \frac{2x}{\lambda} -r\right)\frac{\lambda^k t^k}{k!}
\sum_{m=0}^\infty S_2(m+r,r) 2^{m+r} \frac{1}{\binom{m+r}{r}} \frac{t^m}{m!}  \cr
=&\sum_{n=0}^\infty \left(  \sum_{m=0}^n \binom{n}{m} B_{n-m}^* \left(\frac{2x}{\lambda}-r \right)
\lambda^{n-m} \frac{S_2(m+r,r)}{\binom{m+r}{r}} 2^{m+r}\right)\frac{t^n}{n!}.
\end{align}
On the other hand,
\begin{align}\label{eq20}
\sum_{m=0}^\infty b_{m,\lambda}^{(r)} (x) \frac{1}{m!}(e^{2t}-1)^m
=&\sum_{m=0}^{\infty}b_{m,\lambda}^{(r)} (x)\sum_{n=m}^\infty S_2(n,m) 2^n \frac{t^n}{n!}\cr
=&\sum_{n=0}^{\infty} \left(\sum_{m=0}^{n} b_{m,\lambda}^{(r)}(x)2^n S_2(n,m) \right)\frac{t^n}{n!}.
\end{align}
From \eqref{eq19} and \eqref{eq20}, we have
\begin{equation}\label{eq21}
\sum_{m=0}^n b_{m,\lambda}^{(r)} (x) S_2(n,m)
= \sum_{m=0}^n \binom{n}{m} B_{n-m}^* \left(\frac{2x}{\lambda}-r \right)
\lambda^{n-m} \frac{S_2(m+r,r)}{\binom{m+r}{r}} 2^{m+r-n}.
\end{equation}
Now, we define the {\it{type $2$ degenerate Bernoulli polynomials of the second kind}} by
\begin{equation}\label{eq22}
\frac{(1+t)-(1+t)^{-1}}{\log_\lambda(1+t)} (1+t)^x = \sum_{n=0}^\infty b_{n, \lambda}^* (x) \frac{t^n}{n!}.
\end{equation}
When $x=0$, $b_{n,\lambda}^*=b_{n,\lambda}^*(0)$ are called the type $2$ degenerate Bernoulli numbers of the
second kind.
Note that $\lim_{\lambda\rightarrow 0} b_{n,\lambda}^*(x)=b_n^*(x)$,
where $b_n^*(x)$ are the type $2$ Bernoulli polynomials of the second kind given by 
\begin{equation*}
\frac{(1+t)-(1+t)^{-1}}{\log(1+ t)} (1+ t)^x = \sum_{n=0}^{\infty} b_n^* (x) \frac{t^n}{n!}.
\end{equation*}
From \eqref{eq17} and \eqref{eq22}, we note that
\begin{align}\label{eq23}
\frac{(1+t)-(1+t)^{-1}}{\log_\lambda(1+t)} (1+t)^x &= \frac{t}{\log_{\lambda}(1+t)}(1+t)^x \left( 1+ \frac{1}{1+t} \right)\cr
&= \frac{t}{\log_{\lambda} (1+t)}(1+t)^x + \frac{t}{\log_{\lambda}(1+t)}(1+t)^{x-1}\cr
&= \sum_{n=0}^\infty \left( b_{n,\lambda}^{(1)} (x)+ b_{n,\lambda}^{(1)}(x-1)  \right)\frac{t^n}{n!}.
\end{align}

Therefore, we obtain the following theorem.

\begin{theorem}\label{thm1}
For $n\geq 0$, we have
\begin{equation*}
b_{n,\lambda}^*(x)= b_{n,\lambda}^{(1)} (x)+ b_{n,\lambda}^{(1)}(x-1).
\end{equation*}
Moreover,
\begin{equation*}
\sum_{m=0}^n  b_{m,\lambda}^{(r)} (x)S_2(n,m)
= \sum_{m=0}^n \binom{n}{m} B_{n-m}^*\left( \frac{2x}{\lambda}-r \right)
\lambda^{n-m} \frac{S_2(m+r,r)}{\binom{m+r}{r}} 2^{m+r-n},
\end{equation*}
where $r$ is a positive integer.
\end{theorem}

Now, we observe that
\begin{align}\label{eq24}
\frac{(1+t)-(1+t)^{-1}}{\log_\lambda(1+t)} (1+t)^x
&=\sum_{l=0}^\infty b_{l, \lambda}^*  \frac{t^l}{l!}
\sum_{m=0}^\infty (x)_m \frac{t^m}{m!} \cr
&= \sum_{n=0}^\infty \left( \sum_{l=0}^n \binom{n}{l} b_{l,\lambda}^* (x)_{n-l}  \right)\frac{t^n}{n!},
\end{align}
where $(x)_0=1$, $(x)_n=x(x-1)\cdots (x-n+1)$, $(n\geq 1)$.
From \eqref{eq22} and \eqref{eq24}, we get
\begin{align}\label{eq25}
b_{n,\lambda}^*(x)=\sum_{l=0}^n \binom{n}{l} b_{l,\lambda}^*(x)_{n-l}, \;\;(n\geq 0).
\end{align}

For $\alpha\in\mathbb{R}$, let us define the {\it{type $2$ degenerate Bernoulli polynomials of the second kind of order $\alpha$}} by 
\begin{align}\label{eq26}
\left( \frac{(1+t)-(1+t)^{-1}}{\log_\lambda(1+t)}\right)^\alpha (1+t)^x
=\sum_{n=0}^{\infty} b_{n,\lambda}^{*(\alpha)} (x) \frac{t^n}{n!}
\end{align}
When $x=0$, $b_{n,\lambda}^{*(\alpha)}=b_{n,\lambda}^{*(\alpha)}(0)$ are called the type $2$ degenerate Bernoulli numbers of the
second kind of order $\alpha$.

Let $\alpha=k\in \mathbb{N}$. Then  we have
\begin{align}\label{eq27}
\sum_{n=0}^{\infty} b_{n,\lambda}^{*(k)} (x) \frac{t^n}{n!}
=\left( \frac{(1+t)-(1+t)^{-1}}{\log_\lambda(1+t)}\right)^k (1+t)^x.
\end{align}
By replacing $t$ by $e_\lambda(t)-1$ in \eqref{eq27}, we get
\begin{align}\label{eq28}
\frac{k!}{t^k} \frac{1}{k!} \left(e_\lambda(t)-e_\lambda^{-1}(t)  \right)^k e_\lambda^x(t)
&= \sum_{l=0}^{\infty} b_{l,\lambda}^{*(k)} (x)\frac{1}{l!}(e_\lambda(t)-1)^l \cr
&= \sum_{l=0}^{\infty} b_{l,\lambda}^{*(k)} (x) \sum_{n=l}^\infty S_{2, \lambda}(n,l) \frac{t^n}{n!}\cr
&= \sum_{n=0}^\infty \left(\sum_{l=0}^n b_{l,\lambda}^{*(k)}(x) S_{2,\lambda}(n,l)  \right) \frac{t^n}{n!},
\end{align}
where $S_{2,\lambda}(n,l)$ are the degenerate Stirling numbers of the second kind given by
\begin{align}\label{eq29}
\frac{1}{k!} \left(e_\lambda(t)-1 \right)^k
= \sum_{n=k}^\infty S_{2,\lambda}(n,k) \frac{t^n}{n!}, \;\;\;
\quad \text{ (see \cite{ref10})}.
\end{align}

On the other hand, we also have
\begin{align}\label{eq30}
\frac{k!}{t^k}\frac{1}{k!} \left(e_\lambda(t)-e_\lambda^{-1}(t) \right)^k e_\lambda^x (t)
=&\frac{k!}{t^k}\frac{1}{k!} \left(e_\lambda^2(t)-1 \right)^k e_\lambda^{x-k} (t)    \cr
=&\frac{k!}{t^k}\frac{1}{k!} \left(e_{\frac{\lambda}{2}}(2t)-1 \right)^k e_\lambda^{x-k} (t)    \cr
=& \sum_{m=0}^\infty S_{2,\frac{\lambda}{2}}(m+k,k) \frac{2^{m+k}}{\binom{m+k}{k}} \frac{t^m}{m!} 
 \sum_{l=0}^\infty (x-k)_{l, \lambda} \frac{t^l}{l!} \cr
=& \sum_{n=0}^\infty \left( \sum_{m=0}^n \frac{\binom{n}{m}2^{m+k}}{\binom{m+k}{k}} S_{2, \frac{\lambda}{2}}(m+k,k) (x-k)_{n-m,\lambda} \right)
\frac{t^n}{n!}.
\end{align}
Therefore, by \eqref{eq28} and \eqref{eq30}, we obtain the following theorem.

\begin{theorem}\label{thm2}
For $n\geq 0$, we have
\begin{equation*}
\sum_{l=0}^n  b_{l,\lambda}^{*(k)} (x)S_{2,\lambda}(n,l)
= \sum_{l=0}^n \frac{\binom{n}{l}2^{l+k}}{\binom{l+k}{k}}S_{2,\frac{\lambda}{2}}(l+k,k)  (x-k)_{n-l, \lambda}.
\end{equation*}
\end{theorem}
In particular,
\begin{equation*}
2^{n+k}S_{2,\frac{\lambda}{2}} (n+k,k)
= \binom{n+k}{k} \sum_{l=0}^n b_{l,\lambda}^{*(k)}(k) S_{2,\lambda}(n,l).
\end{equation*}

For $\alpha \in \mathbb{R}$, we recall that the type $2$ degenerate Bernoulli polynomials of
order $\alpha$ are defined by 
\begin{align}\label{eq31}
\left( \frac{t}{e_\lambda(t) -e_\lambda^{-1}(t) }\right)^\alpha e_\lambda^x(t)
= \sum_{n=0}^{\infty}\beta_{n, \lambda}^{*(\alpha)}(x)\frac{t^n}{n!},\quad \text{ (see \cite{ref6, ref12})}.
\end{align}
For $k \in \mathbb{N}$, let us take $\alpha=-k$ and replace $t$ by $\log_\lambda (1+t)$ in \eqref{eq31}.
Then we have
\begin{align}\label{eq32}
\left( \frac{(1+t)-(1+t)^{-1}}{\log_\lambda(1+t)}\right)^k (1+t)^x 
=& \sum_{l=0}^{\infty}\beta_{l, \lambda}^{*(-k)}(x)\frac{1}{l~!}\left(\log_\lambda (1+t) \right)^l\cr
=& \sum_{l=0}^\infty \beta_{l, \lambda}^{*(-k)}(x)\sum_{n=l}^\infty S_{1,\lambda}(n.l) \frac{t^n}{n!}\cr
=& \sum_{n=0}^\infty \left( \sum_{l=0}^n \beta_{l, \lambda}^{*(-k)} S_{1,\lambda}(n.l)\right) \frac{t^n}{n!},
\end{align}
where $S_{1,\lambda}(n,l)$ are the degenerate Stirling numbers of the first kind given by
\begin{equation}\label{eq32-1}
\frac{1}{k!}(\log_{\lambda}(1+t))^k = \sum_{n=k}^{\infty}S_{1,\lambda}(n,k)\frac{t^n}{n!}.
\end{equation}
Note here that 
$\lim_{\lambda\rightarrow 0}S_{1, \lambda}(n,l)=S_1(n,l)$.
Therefore, by \eqref{eq24} and \eqref{eq32}, we obtain the following theorem.

\begin{theorem}\label{thm3}
For $n\geq 0$ and $k \in \mathbb{N}$, we have
\begin{equation*}
b_{n,\lambda}^{*(k)}(x)= \sum_{l=0}^n \beta_{l,\lambda}^{*(-k)} (x) S_{1,\lambda}(n,l).
\end{equation*}
\end{theorem}
We observe that
\begin{align}\label{eq33}
\frac{1}{k!}t^k &= \frac{1}{k!} \left( (1+t)^{\frac{1}{2}} -(1+t)^{-\frac{1}{2}}\right)^k (1+t)^{\frac{k}{2}}\cr
=&  \frac{1}{k!} \left( e_\lambda^{\frac{1}{2}} (\log_\lambda(1+t))-e_\lambda^{-\frac{1}{2}}
\left( \log_\lambda(1+t) \right)\right)^k (1+t)^{\frac{k}{2}}\cr
=& \sum_{l=k}^\infty T_\lambda(l,k) \frac{1}{l!} (\log_\lambda(1+t))^l  
\sum_{r=0}^\infty \left(\frac{k}{2} \right)_r \frac{t^r}{r!}\cr
=& \sum_{l=k}^\infty T_\lambda(l,k)  \sum_{m=l}^{\infty} S_{1,\lambda}(m,l)  \frac{t^m}{m!}       
\sum_{r=0}^\infty \left( \frac{k}{2} \right)_r \frac{t^r}{r!} \cr
=& \sum_{m=k}^\infty \sum_{l=k}^{m} T_\lambda(l,k)  S_{1,\lambda}(m,l)  \frac{t^m}{m!}       
\sum_{r=0}^\infty \left( \frac{k}{2} \right)_r \frac{t^r}{r!} \cr
=& \sum_{n=k}^\infty \left( \sum_{m=k}^n \sum_{l=k}^m T_\lambda(l,k) S_{1,\lambda}(m,l) 
\binom{n}{m} \left(\frac{k}{2} \right)_{n-m}\right)\frac{t^n}{n!}.
\end{align}
On the other hand,
\begin{align}\label{eq34}
\frac{1}{k!}t^k &= \left( \frac{t}{\log_\lambda(1+t)}\right)^k \frac{1}{k!}\left(\log_\lambda(1+t)\right)^k\cr
=& \sum_{l=0}^\infty b_{l,\lambda}^{(k)} \frac{t^l}{l!}  
\sum_{m=k}^\infty S_{1,\lambda}(m,k) \frac{t^m}{m!} \cr
=& \sum_{n=k}^\infty \left( \sum_{m=k}^n S_{1,\lambda}(m,k) b_{n-m,\lambda}^{(k)} \binom{n}{m}  \right)\frac{t^n}{n!}.
\end{align}
Therefore, by \eqref{eq33} and \eqref{eq34}, we obtain the following theorem.

\begin{theorem}\label{thm4}
For $n, k\geq 0$, we have
\begin{equation*}
\sum_{m=k}^n \sum_{l=k}^m T_\lambda(l,k) S_{1,\lambda}(m,l)
\binom{n}{m} \left(\frac{k}{2} \right)_{n-m}=\sum_{m=k}^n S_{1,\lambda}(m,k) b_{n-m,\lambda}^{(k)} \binom{n}{m}.
\end{equation*}
\end{theorem}

\medskip

\section{\bf Conclusions}
\medskip

In [1,2], Carlitz initiated study of the degenerate Bernoulli and Euler polynomials. In recent years, many mathematicians have investigated various degenerate versions of some old and new polynomials and numbers, and found quite a few interesting results [3,4,6,10,11]. It is remarkable that studying degenerate versions is not only limited to polynomials but also can be applied to transcendental functions. Indeed, the degenerate gamma functions were introduced and studied in [8,9]. \\
\indent In this paper, we introduced the type 2 degenerate Bernoulli polynomials of the second kind and their higher-order analogues, and studied some identities and expressions for these polynomials. Specifically, we obtained a relation between the type 2 degenerate Bernoulli polynomials of the second and the degenerate Bernoulli polynomials of the second, an identity involving higher-order analogues of those polynomials and the degenerate Stirling numbers of second kind, and an expression of higher-order analogues of those polynomials in terms of the higher-order type 2 degenerate Bernoulli polynomials and the degenerate Stirling numbers of the first kind. \\
\indent In addition, we obtained an identity involving the higher-order degenerate Bernoulli polynomials of the second kind, the type 2 Bernoulli polynomials and Stirling numbers of the second kind, and 
an identity involving the degenerate central factorial numbers of the second kind, the degenerate Stirling numbers of the first kind and the higher-order degenerate Bernoulli polynomials of the second kind.

\bigskip

{\bf Competing interests:}
The authors declare that they have no competing interests.
\bigskip

{\bf Funding:} This research received no external funding.

\end{document}